\documentclass[11pt]{nyjm}

\title{Generalized improper integral definition for infinite limit} 

\author{Michael A. Blischke}  
\address{} 
\email{mblischke@att.net}   



\newtheorem{thm}{Theorem}

\newtheorem{lem}[thm]{Lemma}

\theoremstyle{definition}
\newtheorem{definition}{Definition}
\theoremstyle{example}
\newtheorem{example}{Example}

\usepackage{url}

\begin{document} 
 
\begin{abstract}  
 A generalization of the definition of a one-dimensional improper integral with an infinite limit is presented.  The new definition extends the range of valid integrals to include integrals which were previously considered to not be integrable.  This definition is shown to preserve linearity and uniqueness.  Integrals which are valid under the conventional definition have the same value under the new definition.  Criteria for interchanging the order of integration and differentiation, and for interchanging the order of integration with a second integration, are determined.  Examples are provided.  A restriction on changing the variable of integration using integration by substitution with the new definition is demonstrated.
\end{abstract} 
\maketitle
\tableofcontents


\section{Introduction and definition}

The standard improper integral defined by
\[
\int_{a}^{\infty} f(x) \,dx \equiv \lim_{b\rightarrow\infty} \left\{ \int_{a}^{b} f(x) \,dx \right\}
\]
exists when the limit exists.  When the limit does not exist, the integral is conventionally said to not exist, or to diverge.  As an example of this case, take $f(x) = \sin(x)$.  Then
\[
\int_{a}^{b} f(x) \,dx = \int_{a}^{b} \sin(x) \,dx = \cos(a) - \cos(b)
\]
which oscillates as $b \to \infty$, and does not approach a single value.

With only the assumption that the integral \emph{does} exist, however, it's easy to show (using $\sin(x+\pi) = -\sin(x)$) that its value would have to be $\cos(a)$.  This suggests that it may be possible to define an integral that applies to improper integrals which oscillate as the limit on the upper bound is taken.

Extending the family of functions which may be integrated is, on its own, sufficient motivation for pursuing this development.  As a more practical motivation, consider a situation where one needs to evaluate
\[
\frac{\partial}{\partial y} \left[\int^{\infty}_{a} f(x,y) \,dx \right]
\]
and would be able to obtain a greatly simplified integrand if the derivative could be brought inside the integral, as
\[
\int^{\infty}_{a} \left[\frac{\partial f(x,y)}{\partial y}\right] \,dx,
\]
but finds that the resulting integral does not exist.

The sharp termination of the finite integral at the artificial upper bound $b$ in the conventional definition is arbitrary, as the desired integral has an unbounded range.  Further, it is that sharp termination that is the source of the oscillations when taking the limit in the example above.  In order to have a single limiting value, the oscillations that occur when taking the limit need to be removed.  

To this end, we propose the following generalization of the standard improper integral definition:
\begin{definition}
\begin{equation}
\label{D01}
\int^{\infty}_{a} \hspace{-19.1 bp} \textsf{Z} \hspace{10 bp} f(x) \,dx \equiv \lim_{b\rightarrow\infty}\left\{\int^{b}_{a}f(x) \,dx + \int^{b+c}_{b} f(x) z(x-b) \,dx\right\}
\end{equation}
\end{definition}
Here $f(x)$ is the function to be integrated, and $z(x)$ is a termination function (defined below).  The inclusion of the additional term inside the limit allows us to rigorously achieve convergence.  To the best of the author's knowledge, the definition in (\ref{D01}) has never been proposed elsewhere.

Throughout this paper, for clarity, we include the over-struck \sffamily Z \normalfont on integrals using the proposed definition and a termination function, to distinguish them from integrals that exist under conventional definitions.  We will assume the Riemann definition for conventional integrals when necessary (including its extension to improper integrals with an infinite limit), but most of the development does not depend on the particular definition.  We will point out where the Riemann definition is assumed.

When the integrals and the limit in (\ref{D01}) exist, the integral of $f(x)$ is said to exist with respect to $z(x)$ under our definition, and to have the value of the limit.  If, for all termination functions for which the limit in (\ref{D01}) exists, the value of that limit is the same, then the integral of $f(x)$ is said to exist, with that value for the integral.

We will examine the case where the upper bound is infinite.  The corresponding case where the lower bound is infinite follows easily, and the case where both bounds are infinite is found as the sum of two integrals with one infinite bound each, in the obvious manner.  We will not consider here the case of improper integrals with finite limits.

We envision the term involving $z(x)$ in (\ref{D01}) as smoothing out the sharp termination of the conventional definition, so $z(x)$ should not be arbitrary.  To that end we require $z(x)$ to satisfy the following:
\begin{equation}
\label{D03}
z(x) = 1 \mbox{ for } x \leq 0,
\end{equation}
\begin{equation}
\label{D04}
z(x) = 0 \normalfont \mbox{ for } x \geq c ,
\end{equation}
and its derivative is required to satisfy
\begin{equation}
\label{D04a}
z'(x) \le 0,
\end{equation}
\begin{equation}
\label{D05}
\int_{-\infty}^{\infty}z'(x) \,dx = \int^{c}_{0} z'(x) \,dx \equiv -1 .
\end{equation}
Conditions (\ref{D03}) through (\ref{D04a}) imply that $z(x)$ is monotonically decreasing, over the range of 1 to 0, for all $x$.

\begin{definition}
A termination function is any function satisfying the conditions given in (\ref{D03}) through (\ref{D05}).  A termination function for $f(x)$ is a termination function for which the limit in (\ref{D01}) exists.  
\end{definition}

Condition (\ref{D05}) follows from (\ref{D03}) through (\ref{D04a}), and so strictly is not required explicitly.  In practice, however, the monotonicity requirement (\ref{D04a}) may be relaxed somewhat to
\begin{equation}
\label{D04b}
0 \le z(x) \le 1 ,
\end{equation}
in which case (\ref{D05}) is required explicitly, hence its inclusion above.  As will be shown later, for any termination function satisfying (\ref{D04b}), a termination function satisfying (\ref{D04a}) may be constructed, so the stricter requirement (\ref{D04a}) is not a burden.  I use the requirement (\ref{D04a}) instead of (\ref{D04b}) here because the following development is easier.

Implicit in the requirements of (\ref{D03}) through (\ref{D05}), is that $z(x)$ is not a function of $b$.  If $z(x)$ were allowed to vary with $b$ in an arbitrary fashion, for some integrals (where conventional integration fails), one could carefully choose $z(x,b)$ to cause the limit in (\ref{D01}) to be any desired value.  Not all forms of variation of termination functions WRT $b$ will necessarily be a problem, and it is possible some restricted forms of variation WRT $b$ may allow us to expand the domain of integrable functions.

The condition (\ref{D04}) may also be more strict than necessary.  A less restrictive condition, perhaps $z(x) \rightarrow 0$ as $x \rightarrow \infty$, coupled with the integral of $z'(x)$ being absolutely convergent, may be sufficient for the developments here.  A looser condition such as that may result in the ability to integrate additional integrals.  

Investigation into less restrictive conditions on termination functions which maintain the desirable properties of the new definition, but which expand the allowable set of integrable functions, is a subject for future research.  For this paper we will proceed using the conditions (\ref{D03}) through (\ref{D05}), and require termination functions have no variation WRT $b$.

The definition in (\ref{D01}) can be manipulated into a more useful form.  Beginning with (\ref{D01}),
\begin{equation}
\int^{\infty}_{a} \hspace{-19.1 bp} \textsf{Z} \hspace{10 bp} f(x) \,dx = \lim_{b\rightarrow\infty}\left\{\left.
\begin{array}{l}
\!\!\!\! \\
\!\!\!\! \end{array}
F(x)\right|^{b}_{a} \,dx + \int^{c}_{0} f(x+b) z(x) \,dx\right\}
\end{equation}
where
\begin{equation}
\label{D07}
F(x) \equiv \int^{x} f(x') \,dx' .
\end{equation}

Using integration by parts,
\begin{equation}
\int^{\infty}_{a} \hspace{-19.1 bp} \textsf{Z} \hspace{9 bp} f(x) \,dx = \lim_{b\to\infty}\left\{\!F(b)\!-\!F(a) + \!\! \left.
\begin{array}{l}
\!\!\!\! \\
\!\!\!\! \end{array}
F(x+b) z(x) \right|^{c}_{0} \!\,dx - \int^{c}_{0} F(x+b) z'(x) \,dx\right\} .
\end{equation}

Now using conditions (\ref{D03}) and (\ref{D04}), we obtain
\begin{equation}
\label{D09}
\int^{\infty}_{a} \hspace{-19.25 bp} \textsf{Z} \hspace{10 bp} f(x) \,dx = -F(a) - \lim_{b\to\infty}\left\{\int^{c}_{0} F(x+b) z'(x) \,dx\right\}
\end{equation}
where $F(a)$ is constant WRT $b$, and is pulled outside the limit.  We will make extensive use of (\ref{D09}).

Given $z(x)$, $z'(x)$ is found by differentiation.  If we are instead given $z'(x)$, $z(x)$ can be found as
\begin{equation}
\label{D10}
z(x) = \left[-z'(x) \right] \otimes \left( 1 - s(x) \right)
\end{equation}
where $\otimes$ denotes convolution, and where $s(x)$ is the unit step function.  Thus, 
\begin{eqnarray}
\label{D11}
z(x) &=& -\int^{\infty}_{-\infty} z'(x') \,dx'  + \int^{\infty}_{-\infty} s(x') z'(x-x') \,dx' \nonumber \\
&=& 1 + \int^{\infty}_{0} z'(x-x') \,dx' \nonumber \\
&=& 1 - \int^{-\infty}_{x} z'(x') \,dx' \nonumber \\
z(x) &=& 1 + \int^{x}_{-\infty} z'(x') \,dx' .
\end{eqnarray}
Taking the derivative of both sides of (\ref{D11}) yields the identity $z'(x) = z'(x)$.

Given two termination functions, a third can be created by convolving their derivatives to get the derivative of the third.  Given termination functions $z_1(x)$ and $z_2(x)$, we have
\begin{equation}
\label{D12}
z_i'(x) = \left\{ \begin{array}{lcc}
0 \mbox{ } \mbox{ } \mbox{ } x < 0 \\
0 \mbox{ } \mbox{ } \mbox{ } x > c_i \\
\end{array} \right.
\end{equation}
with
\begin{equation}
\label{D13}
\int^{\infty}_{-\infty} z_i'(x) \,dx = -1 .
\end{equation}

Define
\begin{equation}
\label{D14}
z'(x) = - z_1'(x) \otimes z_2'(x) \equiv - \int^{\infty}_{-\infty} z_1'(x-x') z_2'(x') \,dx' .
\end{equation}

Using (\ref{D12}), this becomes
\begin{equation}
\label{D15}
z'(x) \equiv - \int^{c_2}_{0} z_1'(x-x') z_2'(x') \,dx' .
\end{equation}
The integrand in (\ref{D15}) is non-zero only when $x-x'$, the argument of $z_1'$, is between $0$ and $c_1$.  Applying also that $z_2'$ is non-zero when $0 \leq x' \leq c_2$, it can be seen that
\begin{equation}
\label{D16}
z'(x) = \left\{ \begin{array}{lcc}
0 \mbox{   } \mbox{ } \mbox{ } x < 0 \\
0 \mbox{   } \mbox{ } \mbox{ } x > c \\
\end{array} \right.
\end{equation}
where $c \equiv c_1 + c_2$.

Since (\ref{D14}) is the convolution integral, using (\ref{D13}) for $z_1$ and $z_2$ immediately gives us
\begin{equation}
\label{D17}
\int^{\infty}_{-\infty} z'(x) \,dx = \int^{c}_{0} z'(x) \,dx \equiv -1
\end{equation}
which satisfies condition (\ref{D05}).  Using (\ref{D16}) in (\ref{D11}), condition (\ref{D03}) is satisfied, and then using (\ref{D17}), we see that condition (\ref{D04}) is also satisfied.  The condition in (\ref{D04a}) is also easily seen to be satisfied.  Thus, $z(x)$ found using (\ref{D14}) and (\ref{D11}) is a termination function.

Note that in place of (\ref{D14}), the equivalent definition
\begin{equation}
z'(x) \equiv - \int^{c_1}_{0} z_1'(x') z_2'(x-x') \,dx'
\end{equation}
could have been used.  We will call $z(x)$ the combined termination function of $z_1(x)$ and $z_2(x)$, and will call $z_1(x)$ and $z_2(x)$ components of the combined termination function $z(x)$.

Using (\ref{D14}) in (\ref{D11}), we can get the combined termination function as
\begin{equation}
z(x) = 1 + \int^{x}_{-\infty} \left[ - \int^{\infty}_{-\infty} z_1'(x'-t) z_2'(t) \,dt \right] \,dx' .
\end{equation}
For convenience, we will write this relation between the combined termination function and its components using the notation
\begin{equation}
\label{D20}
z(x) = z_1(x) \odot z_2(x) .
\end{equation}

\section{Properties of the generalized definition}
In this section, we will develop some important properties of the integral definition.
\begin{thm}
\label{TM01}
For a function $f(x)$, the integral with respect to a termination function $z_1(x)$ gives the same value as the integral with respect to a combined termination function having $z_1(x)$ as one if its components, with an arbitrary termination function $z_2(x)$ as its other component.
\end{thm}

\begin{proof}
Using (\ref{D17}), the limit term from the RHS of (\ref{D09}) becomes
\begin{eqnarray}
& & \hspace{-60 bp} \lim_{b\to\infty}\left\{\int^{c_1}_{0} F(x+b) z_1'(x) \,dx \right\} \nonumber \\
&=& -\left[ \lim_{b\to\infty}\left\{ \int^{c_1}_{0} F(x+b) z_1'(x) \,dx \right\} \right] \int^{c_2}_{0} z_2'(x') \,dx' .
\end{eqnarray}
The limiting operation is invariant WRT $x'$, and can be pulled inside of the integration over $x'$, giving
\begin{equation}
\label{D22}
= -\int^{c_2}_{0} z_2'(x') \left[ \lim_{b\to\infty}\left\{ \int^{c_1}_{0} F(x+b) z_1'(x) \,dx \right\} \right] \,dx' .
\end{equation}

The limit as $b$ approaches $\infty$ is the same as the limit as $x'+b$ approaches $\infty$ for any finite value of $x'$, so (\ref{D22}) can be written as
\begin{equation}
\label{D23}
= -\int^{c_2}_{0} z_2'(x') \left[ \lim_{b\to\infty}\left\{ \int^{c_1}_{0} F(x+x'+b) z_1'(x) \,dx \right\} \right] \,dx' .
\end{equation}
The limit WRT $b$ converges uniformly, so the order of the limit WRT $b$ and the integration WRT $x'$ can be safely interchanged (See for example ~\cite{Spielberg:limits} for this and subsequent interchanges of limiting operations), giving
\begin{equation}
=-\lim_{b\to\infty}\left\{ \int^{c_2}_{0} z_2'(x') \left[ \int^{c_1}_{0} F(x+x'+b) z_1'(x) \,dx \right] \,dx'\right\} .
\end{equation}
Since $x' \geq 0$, and since $z_1'(x) \equiv 0$ for $x \leq 0$, the lower limit of the integration over $x$ can be changed to $-x'$.  Since $z_1'(x) \equiv 0$ for $x \geq c_1$, and $c-x' >= c_1$, the upper limit of the integration over $x$ can be changed to $c-x'$.  Thus
\begin{eqnarray}
&=& -\lim_{b\to\infty}\left\{\int^{c_2}_{0} z_2'(x') \left[ \int^{c-x'}_{-x'} F(x+x'+b) z_1'(x) \,dx \right] \,dx'\right\} \nonumber \\
&=& -\lim_{b\to\infty}\left\{\int^{c_2}_{0} z_2'(x') \left[ \int^{c}_{0} F(x+b) z_1'(x-x') \,dx \right] \,dx'\right\} \nonumber \\
&=& \lim_{b\to\infty}\left\{\int^{c}_{0} F(x+b) \left[ - \int^{c_2}_{0} z_1'(x-x') z_2'(x') \,dx' \right] \,dx \right\} \nonumber \\
&=& \lim_{b\to\infty}\left\{\int^{c}_{0} F(x+b) z'(x) \,dx \right\} .
\end{eqnarray}

We thus obtain, for $z_1(x)$ a termination function of $f(x)$ and a component of the combined termination function $z(x)$,
\begin{eqnarray}
\label{D26}
\lim_{b\to\infty} \left\{ \int^{c_1}_{0} F(x+b) z_1'(x) \,dx \right\} = \lim_{b\to\infty}\left\{\int^{c}_{0} F(x+b) z'(x) \,dx \right\} .
\end{eqnarray}

Using (\ref{D09}), we then have that
\begin{eqnarray}
\int^{\infty}_{a} \hspace{-19.1 bp} \textsf{Z} f(x) \,dx &=& -F(a) - \left[ \lim_{b\to\infty} \left\{ \int^{c_1}_{0} f(x+b) z_1'(x) \,dx \right\} \right] \nonumber \\
&=& -F(a) - \left[ \lim_{b\to\infty} \left\{ \int^{c}_{0} f(x+b) z'(x) \,dx \right\} \right] .
\end{eqnarray}
\end{proof}
Theorem \ref{TM01} is crucial result.  Given a termination function $z_1(x)$ for a function $f(x)$, and given \textit{any} other termination function $z_2(x)$, their combined termination function $z(x)$ is also a termination function for $f(x)$, and the integral of $f(x)$ using either of those two termination functions \textit{has the same value}.  

Given any two different termination functions $z_1(x)$ and $z_2(x)$ for $f(x)$ (i.e. for which the integral (\ref{D01}) exists) their combined termination function is also a termination function for $f(x)$, and the resulting value of the integration is equal to that of $z_1(x)$.  Since the choice of which termination function is $z_1(x)$ and which is $z_2(x)$ is arbitrary, a function which is integrable WRT some termination function $z(x)$ has the same value for \textit{all} termination functions for which the limit exists.  Thus, if $f(x)$ is integrable under our definition WRT any termination function $z(x)$, we can simply say that $f(x)$ is integrable under our definition, without further qualification.  That is, our definition gives us a unique answer.

We immediately see that conventional integration with an infinite upper bound is a special case of the general definition:
\begin{lem}
When the improper integral exists using the conventional improper integral definition,
\begin{equation}
\int^{\infty}_{a} \hspace{-19.8 bp} \textsf{Z} \hspace{10 bp} f(x) \,dx = \int^{\infty}_{a} f(x) \,dx .
\end{equation}
\end{lem}
\begin{proof}
Take as the termination function
\begin{equation}
z(x) = \left\{ \begin{array}{lcc}
1 \mbox{   } \mbox{ } \mbox{ } x \leq 0 \\
0 \mbox{   } \mbox{ } \mbox{ } x > 0 \\
\end{array} \right.
\end{equation}
Then
\begin{eqnarray}
\int^{\infty}_{a} \hspace{-19.1 bp} \textsf{Z} \hspace{10 bp} f(x) \,dx &=& \lim_{b\rightarrow\infty}\left\{\int^{b}_{a}f(x) \,dx + \int^{b+c}_{b} f(x) z(x-b) \,dx\right\} \nonumber \\
&=&\lim_{b\rightarrow\infty}\left\{\int^{b}_{a}f(x) \,dx \right\} \nonumber \\
&=&\int^{\infty}_{a} f(x) \,dx .
\end{eqnarray}
\end{proof}
Thus our generalized definition of integration gives the same result as conventional integration, whenever the conventional integral exists.

Additionally, we obtain that, given any two functions $f(x)$ and $g(x)$ integrable under our definition, a linear combination of them is also integrable.  
\begin{thm}
The integral definition (\ref{D01}) satisfies linearity:
\begin{equation}
\label{D36}
\alpha \int^{\infty}_{a} \hspace{-19.9 bp} \textsf{Z} f(x) \,dx + \beta \int^{\infty}_{a} \hspace{-19.9 bp} \textsf{Z} g(x) \,dx \\
\hspace{2 bp} = \int^{\infty}_{a} \hspace{-19.9 bp} \textsf{Z} \left[ \alpha f(x) + \beta g(x) \right] \,dx .
\end{equation}
\end{thm}

\begin{proof}
Assume both integrals on the LHS of (\ref{D36}) exist.
From (\ref{D09})
\begin{equation}
\int^{\infty}_{a} \hspace{-19.1 bp} \textsf{Z} {f(x) \,dx} \\
= -F(a) + \left[ \lim_{b\to\infty}\left\{ \int^{c_f}_{0} F(x+b) z_f'(x) \,dx \right\} \right]
\end{equation}
and
\begin{equation}
\int^{\infty}_{a} \hspace{-19.1 bp} \textsf{Z} g(x) \,dx \\
= -G(a) + \left[ \lim_{b\to\infty}\left\{ \int^{c_g}_{0} G(x+b) z_g'(x) \,dx \right\} \right]
\end{equation}

\begin{eqnarray}
\alpha \int^{\infty}_{a} \hspace{-19.1 bp} \textsf{Z} f(x) \,dx&+&\beta \int^{\infty}_{a} \hspace{-19.1 bp} \textsf{Z} g(x) \,dx = \nonumber \\
& & \alpha \left[-F(a) - \left[ \lim_{b\to\infty}\left\{ \int^{c_f}_{0} F(x+b) z_f'(x) \,dx \right\} \right] \right] \nonumber \\
&+& \beta \left[-G(a) - \left[ \lim_{b\to\infty}\left\{ \int^{c_g}_{0} G(x+b) z_g'(x) \,dx \right\} \right] \right]
\end{eqnarray}
or
\begin{eqnarray}
\label{D33}
\alpha \int^{\infty}_{a} \hspace{-19.1 bp} \textsf{Z} f(x) \,dx + \beta \int^{\infty}_{a} \hspace{-19.1 bp} \textsf{Z} g(x) \,dx &=& -\left[\alpha F(a) + \beta G(a) \right] \nonumber \\
& & \hspace{-110 bp} -\lim_{b\to\infty} \left\{ \alpha \int^{c_f}_{0} F(x+b) z_f'(x) \,dx + \beta \int^{c_g}_{0} G(x+b) z_g'(x) \,dx \right\} .
\end{eqnarray}
Using (\ref{D14}), define
\begin{equation}
z'(x) \equiv - \int^{\infty}_{-\infty} z_f'(x-x') z_g(x') \,dx' ,
\end{equation}
then using (\ref{D26}) twice (once with $z_f'(x)$ = $z_1'(x)$, and once with $z_g'(x)$ = $z_1'(x)$), (\ref{D33}) becomes
\begin{eqnarray}
\alpha \int^{\infty}_{a} \hspace{-19.1 bp} \textsf{Z} f(x) \,dx + \beta \int^{\infty}_{a} \hspace{-19.1 bp} \textsf{Z} g(x) \,dx \hspace{2 bp} &=& -\left[\alpha F(a) + \beta G(a) \right] \nonumber \\
& & \hspace{-70 bp} -\lim_{b\to\infty} \left\{ \int^{c_f+c_g}_{0} \left[ \alpha F(x+b) + \beta G(x+b) \right] z'(x) \,dx \right\}
\end{eqnarray}
Using (\ref{D09}) on the RHS, we then obtain
\begin{equation}
\alpha \int^{\infty}_{a} \hspace{-19.1 bp} \textsf{Z} f(x) \,dx + \beta \int^{\infty}_{a} \hspace{-19.1 bp} \textsf{Z} g(x) \,dx \hspace{2 bp} = \int^{\infty}_{a} \hspace{-19.1 bp} \textsf{Z} \left[ \alpha f(x) + \beta g(x) \right] \,dx
\end{equation}
\end{proof}

Finally, recall our discussion regarding the alternative requirement (\ref{D04b}).  If we have a termination function $z_1'(x)$ satisfying (\ref{D04b}) instead of (\ref{D04a}), a termination function that is the combined termination function of $z_1(x)$ and $z_2(x)$, with $z_2(x)$ defined by 
\begin{equation}
z_2(x) = 1 - \frac{x}{c_1}  \mbox{ } \mbox{ } \mbox{ } \mbox{ } \mbox{ } \mbox{ }  0 \le x \le c_1,
\end{equation}
will satisfy (\ref{D04a}).  Requirements (\ref{D04a}) and (\ref{D04b}) are functionally equivalent, and allow the same integrations to be performed.

\subsection{Examples}

The following shows some example evaluations of integrals which do not exist under the conventional definition.  If conventional improper integral definition were used, example 1 would have constant variation as the upper limit of the integral approached infinity.  Examples \ref{EX02} and \ref{EX03} would have unbounded variation.

\begin{example}
\label{EX01}
\[
f(x) = \sin(\alpha x)
\]
\[
F(x) = \frac{-\cos(\alpha x)}{\alpha}
\]
\[
z(x) = \left\{ \begin{array}{lcc}
1/2 \mbox{   } \mbox{ } \mbox{ } 0 < x < \pi / \alpha \\
0 \mbox{   } \mbox{ } \mbox{ } \mbox{ } \mbox{ } x \geq \pi / \alpha \\
\end{array} \right.
\]
\[
z'(x) = -1/2 \left( \delta(x) + \delta(x-\pi / \alpha) \right)
\]
\begin{eqnarray*}
\int^{\infty}_{a} \hspace{-19.8 bp} \textsf{Z} \hspace{10 bp} \sin(\alpha x) \,dx &=& \frac{\cos(\alpha a)}{\alpha} \\
& & -\frac{1}{2} \lim_{b\rightarrow\infty}\left\{\int_{0}^{\pi}\frac{\cos(\alpha(x+b))}{\alpha} \left( \delta(x) + \delta(x-\pi / \alpha) \right) \,dx \right\} \\
&=& \frac{\cos(\alpha a)}{\alpha} -1/2 \lim_{b\rightarrow\infty}\left\{\frac{\cos(\alpha b) + \cos(\alpha b + \pi)}{\alpha} \right\} \\
&=& \frac{\cos(\alpha a)}{\alpha}
\end{eqnarray*}
\end{example}

\begin{example}
\label{EX02}
\[
f(x) = x \cos(\alpha x)
\]
\[
F(x) = \frac{\cos(\alpha x)}{\alpha^2} + \frac{x\sin(\alpha x)}{\alpha}
\]
\[
z(x) = \left\{ \begin{array}{lcc}
3/4 \mbox{   } \mbox{ } \mbox{ } 0 < x < \pi / \alpha \\
1/4 \mbox{   } \mbox{ } \mbox{ } \pi / \alpha < x < 2\pi / \alpha \\
0 \mbox{   } \mbox{ } \mbox{ } \mbox{ } \mbox{ } x \geq 2\pi / \alpha \\
\end{array} \right.
\]
\[
z'(x) = -1/4 \left( \delta(x) + 2\delta(x-\pi / \alpha) + \delta(x-2\pi / \alpha) \right)
\]
\begin{eqnarray*}
\int^{\infty}_{a} \hspace{-19.8 bp} \textsf{Z} \hspace{10 bp} x\cos(\alpha x) \,dx &=& -\frac{\cos(\alpha a)}{\alpha^2} - \frac{a \sin(\alpha a)}{\alpha} - \\
& &  \hspace{-60 bp} \lim_{b\rightarrow\infty}\left\{\frac{-1}{4}\int_{0}^{\pi/\alpha}\left(\frac{\cos(\alpha(x+b))}{\alpha^2} + \frac{(x+b)\sin(\alpha (x+b))}{\alpha} \right) \right. \cdot \\
& & \left. \begin{array}{l}
\!\!\!\! \\
\!\!\!\! \end{array}
 \hspace{20 bp} \left( \delta(x) + 2\delta(x-\pi / \alpha) + \delta(x-2\pi / \alpha) \right) \,dx \right\}
\end{eqnarray*}
\begin{eqnarray*}
& & \hspace{10 bp} = -\frac{\cos(\alpha a)}{\alpha^2} - \frac{a \sin(\alpha a)}{\alpha} + \\
& & \hspace{20 bp} 1/4 \lim_{b\rightarrow\infty}\left\{ \frac{\cos(\alpha b)}{\alpha^2} + \frac{b\sin (\alpha b)}{\alpha} \right. + \\
& & \hspace{70 bp} \left. 2\left(\frac{\cos(\alpha b+\pi)}{\alpha^2} + \frac{(\pi / \alpha + b)\sin (\alpha b+\pi)}{\alpha}\right) + \right. \\
& & \hspace{70 bp} \left. \frac{\cos(\alpha b+2\pi)}{\alpha^2} + \frac{(2\pi / \alpha + b)\sin (\alpha b+2\pi)}{\alpha} \right\}
\end{eqnarray*}

The trigonometric terms within the limit WRT $b$ cancel, leaving
\begin{eqnarray*}
\int^{\infty}_{a} \hspace{-19.8 bp} \textsf{Z} \hspace{10 bp} x\cos(\alpha x) \,dx &=& -\frac{\cos(\alpha a)}{\alpha^2} - \frac{a \sin(\alpha a)}{\alpha} .
\end{eqnarray*}
\end{example}

\begin{example}
\label{EX03}
\[
f(x) = e^{\beta x}\sin(\alpha x)
\]
\[
F(x) = \frac{e^{\beta x}}{\alpha^2+\beta^2} \left(\beta \sin(\alpha x) -\alpha \cos(\alpha x) \right)
\]
\[
z(x) = \left\{ \begin{array}{lcc}
\frac{1}{1+e^{\beta \pi / \alpha}} \mbox{   } \mbox{ } \mbox{ } 0 < x < \pi / \alpha \\
0 \mbox{   } \mbox{ } \mbox{ } \mbox{ } \mbox{ } x \geq \pi / \alpha \\
\end{array} \right.
\]
\[
z'(x) = -\frac{1}{1+e^{\beta \pi / \alpha}} \left(e^{\beta \pi / \alpha} \delta(x) + \delta(x-\pi / \alpha) \right)
\]
\begin{eqnarray*}
\int^{\infty}_{a} \hspace{-19.8 bp} \textsf{Z} \hspace{10 bp} e^{\beta x}\sin(\alpha x) \,dx &=& -\frac{e^{\beta a}}{\alpha^2+\beta^2} \left(\beta \sin(\alpha a) -\alpha \cos(\alpha a) \right) - \\
& & \hspace{-80 bp} \frac{1}{1+e^{\beta \pi / \alpha}} \lim_{b\rightarrow\infty}\left\{\int_{0}^{\pi / \alpha}\frac{e^{\beta (x+b)}}{\alpha^2+\beta^2} \left(\beta \sin(\alpha (x+b)) -\alpha \cos(\alpha (x+b)) \right) \right.  \cdot \\
& & \left. \hspace{30 bp} \frac{1}{1+e^{\beta \pi / \alpha}} \left(e^{\beta \pi / \alpha} \delta(x) + \delta(x-\pi / \alpha) \right) \,dx \right\}
\end{eqnarray*}

Working through the algebra, all the terms which vary with $b$ cancel completely, leaving
\begin{eqnarray*}
\int^{\infty}_{a} \hspace{-19.8 bp} \textsf{Z} \hspace{10 bp} e^{\beta x}\sin(\alpha x) \,dx&=& -\frac{e^{\beta a}}{\alpha^2+\beta^2} \left(\beta \sin(\alpha a) -\alpha \cos(\alpha a) \right) .
\end{eqnarray*}

When $\beta = 0$, the result is identical to the result from Example \ref{EX01}.  Also, when $\beta < 0$, the result matches the result from conventional integration.
\end{example}

\section{Interchange of order of integration and differentiation}

We will now examine Leibniz's integral rule for differentiation under the integral sign.  Under conventional integration, when certain conditions are met, we have that 
\begin{equation}
\label{D37}
\frac{d}{dy} \left[\int^{\infty}_{a} f(x,y) \,dx \right] = \int^{\infty}_{a} \left[\frac{\partial f(x,y)}{\partial y}\right] \,dx .
\end{equation}
We will show that this relation also holds when one or both integrals require our generalized definition.

Assume
\begin{equation}
\label{D38}
\int^{\infty}_{a} \hspace{-19.1 bp} \textsf{Z} \hspace{10 bp} f(x,y) \,dx = -F(a,y) - \lim_{b\to\infty}\left\{\int^{c_x(y)}_{0} F(x+b,y) z_x(x,y) \,dx\right\}
\end{equation}
exists, with $f(x,y)$ finite and continuous WRT $y$, over a small neighborhood of $y$.  We will also assume that we have a termination function that satisfies\footnote{While at first this may seem a severe restriction, we can often use (\ref{D14}) to obtain a termination function that satisfies (\ref{D38b}) from one which doesn't.  See Example \ref{EX05}.} 
\begin{equation}
\label{D38b}
z_x(x,y) \mbox{ continuous WRT } y .
\end{equation}
Finally, assume
\begin{equation}
\label{D39}
\int^{\infty}_{a} \hspace{-19.1 bp} \textsf{Z} \hspace{10 bp} \left[\frac{\partial f(x,y)}{\partial y}\right] \,dx = -F_y(a,y) - \lim_{b\to\infty}\left\{\int^{\widetilde{c}_x(y)}_{0} F_y(x+b,y) \widetilde{z}_x(x,y) \,dx\right\}
\end{equation}
exists, where
\begin{equation}
\label{D39b}
F_y(x,y) = \frac{\partial F(x,y)}{\partial y} .
\end{equation}

We want to show that (\ref{D37}) holds for the new definition.
\begin{thm}
Assume (\ref{D38}) through (\ref{D39b}) hold.  Then
\begin{equation}
\label{D41}
\frac{d}{dy} \left[\int^{\infty}_{a} \hspace{-19.9 bp} \textsf{Z} \hspace{10 bp} f(x,y) \,dx \right] = \int^{\infty}_{a} \hspace{-19.9 bp} \textsf{Z} \hspace{10 bp} \left[\frac{\partial f(x,y)}{\partial y}\right] \,dx .
\end{equation}
\end{thm}

\begin{proof}
We can always choose a constant value $c$ larger than the maximum of the values $c_x(y)$ needed in (\ref{D37}) and $\widetilde{c}_x(y)$ needed in (\ref{D38}).  For convenience, in the following we will simply use $c$ for the integration limit, and will only alter the upper limit based on how many termination functions are being combined.  We will make extensive use of (\ref{D26}).

First, looking at the LHS of (\ref{D41}),
\begin{equation}
\frac{d}{dy} \left[\int^{\infty}_{a} \hspace{-19.1 bp} \textsf{Z} \hspace{10 bp} f(x,y) \,dx \right] = \lim_{h\to0} \left[ \frac{\int^{\infty}_{a} \hspace{-16.7 bp} \textsf{z} \hspace{10 bp} f(x,y+h) \,dx - \int^{\infty}_{a} \hspace{-16.7 bp} \textsf{z} \hspace{10 bp} f(x,y) \,dx}{h} \right]
\end{equation}
\begin{eqnarray}
&=& \lim_{h\to0} \frac{1}{h}\left[ -F(a,y+h) - \lim_{b\to\infty}\left\{\int^{c}_{0} F(x+b,y+h) z_x(x,y+h) \,dx\right\} + \right. \nonumber \\
& & \hspace{35 bp} \left. F(a,y) + \lim_{b\to\infty}\left\{\int^{c}_{0} F(x+b,y) z_x(x,y) \,dx\right\} \right] \\
&=& \lim_{h\to0} \frac{1}{h}\left[ -F(a,y+h) + F(a,y) - 
\begin{array}{l}
\!\!\!\! \\
\!\!\!\! \end{array}
\right. \nonumber \\
& & \hspace{11 bp} \lim_{b\to\infty}\left\{\int^{3c}_{0} F(x+b,y+h) z_x(x,y+h) \odot z_x(x,y) \odot \widetilde{z}_x(x,y) \,dx\right\} + \nonumber \\
& & \hspace{11 bp} \left. \lim_{b\to\infty}\left\{\int^{3c}_{0} F(x+b,y) z_x(x,y+h) \odot z_x(x,y) \odot \widetilde{z}_x(x,y) \,dx\right\} 
\begin{array}{l}
\!\!\!\! \\
\!\!\!\! \end{array}
\right] .
\end{eqnarray}
Combining the limits on $b$,
\begin{eqnarray}
&=& \lim_{h\to0} \left[ \frac{F(a,y)-F(a,y+h)}{h} - \right. \nonumber \\
& & \left. \frac{1}{h}\lim_{b\to\infty}\left\{\int^{3c}_{0} \left(F(x+b,y+h)-F(x+b,y)\right) \right. \right. \nonumber \\
& & \hspace{70 bp} \left. \left. 
\begin{array}{l}
\!\!\!\! \\
\!\!\!\! \end{array}
z_x(x,y+h) \odot z_x(x,y) \odot \widetilde{z}_x(x,y) \,dx\right\} \right]
\end{eqnarray}
\begin{eqnarray}
&=& -F_y(a,y) - \nonumber \\
& & \lim_{h\to0} \left[ \lim_{b\to\infty}\left\{\int^{3c}_{0} \frac{\left(F(x+b,y+h)-F(x+b,y)\right)}{h} \right. \right. \nonumber \\
& & \hspace{70 bp} \left. \left. 
\begin{array}{l}
\!\!\!\! \\
\!\!\!\! \end{array}
z_x(x,y+h) \odot z_x(x,y) \odot \widetilde{z}_x(x,y) \,dx\right\} \right] .
\end{eqnarray}
The limit WRT $b$ converges uniformly, and the limit WRT $h$ exists, so we can interchange the order of the limits,
\begin{eqnarray}
&=& -F_y(a,y) - \nonumber \\
& & \lim_{b\to\infty}\left\{\lim_{h\to0} \left[ \int^{3c}_{0} \frac{\left(F(x+b,y+h)-F(x+b,y)\right)}{h} \right. \right. \nonumber \\
& & \hspace{70 bp} \left. \left. 
\begin{array}{l}
\!\!\!\! \\
\!\!\!\! \end{array}
z_x(x,y+h) \odot z_x(x,y) \odot \widetilde{z}_x(x,y) \,dx \right]\right\}
\end{eqnarray}

\begin{eqnarray}
&=& -F_y(a,y) - \nonumber \\
& & \lim_{b\to\infty}\left\{ \int^{3c}_{0} \lim_{h\to0} \left[\frac{\left(F(x+b,y+h)-F(x+b,y)\right)}{h} \right. \right. \nonumber \\
& & \hspace{70 bp} \left. \left. 
\begin{array}{l}
\!\!\!\! \\
\!\!\!\! \end{array}
z_x(x,y+h) \odot z_x(x,y) \odot \widetilde{z}_x(x,y) \right] \,dx\right\} ,
\end{eqnarray}
or finally
\begin{eqnarray}
\label{D49}
& & \hspace{-30 bp} \frac{d}{dy} \left[\int^{\infty}_{a} \hspace{-19.1 bp} \textsf{Z} \hspace{10 bp} f(x,y) \,dx \right] = -F_y(a,y) - \nonumber \\
& &  \hspace{-1 bp} \lim_{b\to\infty}\left\{ \int^{3c}_{0} F_y(x+b,y) z_x(x,y) \odot z_x(x,y) \odot \widetilde{z}_x(x,y) \,dx\right\} .
\end{eqnarray}

Now looking at the RHS of (\ref{D41}), we have 
\begin{equation}
\int^{\infty}_{a} \hspace{-19.1 bp} \textsf{Z} \hspace{10 bp} \left[\frac{\partial f(x,y)}{\partial y}\right] \,dx = -F_y(x,y) -\lim_{b\to\infty}\left\{\int^{3c}_{0} F_y(x+b,y) \widetilde{z}_x(x,y) \,dx \right\}
\end{equation}
or
\begin{eqnarray}
\label{D51}
& & \hspace{-30 bp} \int^{\infty}_{a} \hspace{-19.1 bp} \textsf{Z} \hspace{10 bp} \left[\frac{\partial f(x,y)}{\partial y}\right] \,dx = -F_y(x,y) - \nonumber \\
& & \lim_{b\to\infty}\left\{\int^{3c}_{0} F_y(x+b,y) z_x(x,y) \odot z_x(x,y) \odot \widetilde{z}_x(x,y) \,dx \right\} .
\end{eqnarray}

Noting that (\ref{D49}) and (\ref{D51}) are identical, we obtain (\ref{D41}).
\end{proof}
It is an open question whether the existence of the LHS of (\ref{D41}) implies the existence of the RHS.  An affirmative answer would be an important development.

Similarly, while the existence of the integral on the RHS of (\ref{D41}) does not imply the existence of the integral on the LHS (consider $f(x,y)$ equal to a constant), it is perhaps possible that there is always some function $g(x)$ such that
\begin{eqnarray}
\label{D53b}
\int^{\infty}_{a} \hspace{-19.1 bp} \textsf{Z} \hspace{10 bp} \left[\frac{\partial f(x,y)}{\partial y}\right] \,dx &=& \int^{\infty}_{a} \hspace{-19.1 bp} \textsf{Z} \hspace{10 bp} \left[\frac{\partial \left( f(x,y)+g(x) \right) }{\partial y}\right] \,dx \nonumber \\
&=& \frac{d}{dy} \left[\int^{\infty}_{a} \hspace{-19.1 bp} \textsf{Z} \hspace{10 bp} \left(f(x,y)+g(x)\right) \,dx \right].
\end{eqnarray}
Again, this is an open question.

Since conventional integration with an infinite upper bound is a special case of integration using the generalized definition, we also immediately obtain
\begin{lem}
When the corresponding conventional and generalized integrations exist,
\begin{equation}
\label{D52}
\frac{d}{dy} \left[\int^{\infty}_{a} f(x,y) \,dx \right] = \int^{\infty}_{a} \hspace{-19.9 bp} \textsf{Z} \hspace{10 bp} \left[\frac{\partial f(x,y)}{\partial y}\right] \,dx
\end{equation}
\end{lem}
This allows derivative operators to be ``brought inside'' an integral operator in cases where it is not allowed under conventional integration.  This is an important result of this definition of integration.

We also obtain
\begin{lem}
When the corresponding conventional and generalized integrations exist,
\begin{equation}
\frac{d}{dy} \left[\int^{\infty}_{a} \hspace{-19.9 bp} \textsf{Z} \hspace{10 bp} f(x,y) \,dx \right] = \int^{\infty}_{a} \left[\frac{\partial f(x,y)}{\partial y}\right] \,dx .
\end{equation}
\end{lem}

\subsection{Examples involving differentiation}

The following examples show passing differentiation inside the integration.  In the first example, one of the integrals can be performed conventionally, while in the second, both integrations require the general definition.

\begin{example}
\label{EX04}
\[
f(x,y) = \frac{\cos(x y)}{x}
\]
\[
f_y(x,y) = \frac{\partial}{\partial y}\left(\frac{\cos(x y)}{x}\right) = -\sin(xy)
\]
\[
F_y(x,y) = \frac{\cos(x y)}{y}
\]
\[
\widetilde{z}(x) = \left\{ \begin{array}{lcc}
1/2 \mbox{   } \mbox{ } \mbox{ } 0 < x < \pi / y \\
0 \mbox{   } \mbox{ } \mbox{ } \mbox{ } \mbox{ } x \geq \pi / y \\
\end{array} \right.
\]
\[
\widetilde{z}_x(x) = -1/2 \left( \delta(x) + \delta(x-\pi / y) \right)
\]
\begin{eqnarray*}
\int^{\infty}_{a} f(x,y) \,dx = \int^{\infty}_{a} \frac{\cos(x y)}{x} \,dx = -Ci(ay)
\end{eqnarray*}
where ~\cite{abramowitz:cosint}
\[
Ci(z) \equiv \int^{\infty}_{z} \frac{\cos(t)}{t} \,dt
\]

The derivative of $Ci(x)$ is given by~\cite{wolfram:cosint}, Eq. (6)
\begin{eqnarray*}
\frac{d Ci(u)}{du} = -\frac{\cos(u)}{u}
\end{eqnarray*}
so
\begin{eqnarray*}
\frac{d}{dy}\int^{\infty}_{a} \frac{\cos(x y)}{x} \,dx = -\frac{d Ci(ay)}{dy} = -\frac{\cos(a y)}{y}
\end{eqnarray*}
Using Example \ref{EX01}, 
\begin{eqnarray*}
\int^{\infty}_{a} \hspace{-19.8 bp} \textsf{Z} \hspace{10 bp} \frac{\partial}{\partial y}\left(\frac{\cos(x y)}{x}\right) \,dx &=&\int^{\infty}_{a} \hspace{-19.8 bp} \textsf{Z} \hspace{10 bp} -\sin(x y) \,dx \\
&=&-\frac{\cos(ay)}{y}
\end{eqnarray*}
so 
\begin{eqnarray*}
\frac{d}{dy}\int^{\infty}_{a} \frac{\cos(x y)}{x} \,dx = \int^{\infty}_{a} \hspace{-19.8 bp} \textsf{Z} \hspace{10 bp} \frac{\partial}{\partial y}\left(\frac{\cos(x y)}{x}\right) \,dx
\end{eqnarray*}
and we see that (\ref{D52}) is satisfied.
\end{example}

\begin{example}
\label{EX05}
\[
f(x,y) = \sin(x y)
\]
\[
F(x,y) = \frac{-\cos(x y)}{y}
\]
\[
f_y(x,y) = x\cos(x y)
\]
\[
F_y(x,y) = \frac{\cos(x y)}{y^2} + \frac{x\sin(x y)}{x}
\]
\[
z(x,y) = \left\{ \begin{array}{lcc}
1/2 \mbox{   } \mbox{ } \mbox{ } 0 < x < \pi / y \\
0 \mbox{   } \mbox{ } \mbox{ } \mbox{ } \mbox{ } x \geq \pi / y \\
\end{array} \right.
\]
\[
z_x(x,y) = -1/2 \left( \delta(x) + \delta(x-\pi / y) \right) \mbox{ } \footnote{To satisfy (\ref{D38b}), we could combine $z(x,y)$ with $r(x,y)$ and $s(x,y)$, defined by
\[
r_x(x,y) = \left\{ \begin{array}{lcc}
{-y / \pi} {\mbox{   } \mbox{ } \mbox{ }} {0 < x < \pi / y} \\
0 {\mbox{   } \mbox{ } \mbox{ } \mbox{ } \mbox{ }} \mbox{ else} \\
\end{array} \right.
\]
\[
s_x(x,y) = \left\{ \begin{array}{lcc}
{-(3*y_{\max}-2y) / \pi} {\mbox{   } \mbox{ } \mbox{ }} {0 < x < \pi / (3*y_{\max}-2y)} \\
0 {\mbox{   } \mbox{ } \mbox{ } \mbox{ } \mbox{ }} \mbox{ else} \\
\end{array} \right.
\]
where $y_{\max}$ is the maximum $y$ in the domain of interest.  The combination with $r(x,y)$ gets us a termination function whose derivative is a rectangular pulse of width $2y/\pi$.  Combining with $s(x,y)$ then gets us a termination function whose derivative WRT $x$ is non-zero over a constant domain of $3*y_{\max} / \pi$, and continuous WRT $y$, satisfying (\ref{D38b}).}
\]
\[
\widetilde{z}(x,y) = \left\{ \begin{array}{lcc}
3/4 \mbox{   } \mbox{ } \mbox{ } 0 < x < \pi / y \\
1/4 \mbox{   } \mbox{ } \mbox{ } \pi / y < x < 2\pi / y \\
0 \mbox{   } \mbox{ } \mbox{ } \mbox{ } \mbox{ } x \geq 2\pi / y \\
\end{array} \right.
\]
\[
\widetilde{z}_x(x,y) = -1/4 \left( \delta(x) + 2\delta(x-\pi / y) + \delta(x-2\pi / y) \right)
\]

From Example \ref{EX01},
\begin{eqnarray*}
\int^{\infty}_{a} \hspace{-19.8 bp} \textsf{Z} \hspace{10 bp} \sin(x y) \,dx = \frac{\cos(ay)}{y}
\end{eqnarray*}
so
\begin{eqnarray*}
\frac{d}{dy} \left[ \int^{\infty}_{a} \hspace{-19.8 bp} \textsf{Z} \hspace{10 bp} \sin(x y) \,dx \right] &=& \frac{d}{dy} \left[ \frac{\cos(ay)}{y} \right] \\
&=& -\frac{\cos(a y)}{y^2} - \frac{a\sin(a y)}{x}
\end{eqnarray*}

Using Example \ref{EX02}, 
\begin{eqnarray*}
\int^{\infty}_{a} \hspace{-19.8 bp} \textsf{Z} \hspace{10 bp} \frac{\partial \sin(xy)}{\partial y} \,dx &=& 
\int^{\infty}_{a} \hspace{-19.8 bp} \textsf{Z} \hspace{10 bp} x\cos(xy) \,dx \\
&=& -\frac{\cos(ay)}{y^2} - \frac{a \sin(ay)}{y}
\end{eqnarray*}
comparing, we see that
\[
\frac{d}{dy} \left[ \int^{\infty}_{a} \hspace{-19.8 bp} \textsf{Z} \hspace{10 bp} \sin(x y) \,dx \right] = \int^{\infty}_{a} \hspace{-19.8 bp} \textsf{Z} \hspace{10 bp} \frac{\partial \sin(xy)}{\partial y} \,dx,
\]
and (\ref{D41}) is satisfied.
\end{example}

\section{Interchange of order of integration}

We will now examine iterated integration of functions of two variables with respect to both of the variables, where one of the integrals requires our new definition.  The development of this section explicitly assumes Riemann integration for the other integration, and would need to be modified to remove that dependence.

Assume
\begin{equation}
\label{D54}
\int^{\infty}_{a} \hspace{-19.1 bp} \textsf{Z} \hspace{10 bp} f(x,y) \,dx
\end{equation}
exists over the domain $\alpha \leq y \leq \beta$ for some termination function $z(x,y)$.  Assume also that $w(y)$ exists over the same domain, and that
\begin{equation}
\label{D55}
g(x,t) \equiv \int^{t} f(x,y) w(y) \,dy
\end{equation}
exists over the domain $\alpha \leq t \leq \beta$, for $a \leq x < \infty$.  Further, assume that
\begin{equation}
\label{D56}
\int^{\infty}_{a} \hspace{-19.1 bp} \textsf{Z} \hspace{10 bp} g(x,y) \,dx
\end{equation}
exists over the domain $\alpha \leq y \leq \beta$ for some termination function $\widetilde{z}(x,y)$.  It may be that (\ref{D54}) and (\ref{D55}) imply the existence of (\ref{D56}), but that is not proven here.

Analogous to (\ref{D07}), define
\begin{equation}
\label{D58}
F(x,y) \equiv \int^{x} f(x',y) \,dx'
\end{equation}
and
\begin{equation}
\label{D59}
G(x,y) \equiv \int^{x} g(x',y) \,dx' .
\end{equation}

We will also limit ourselves to functions $f(x,y)$ such that
\begin{equation}
\label{D60}
G(x,t) = \int^{t} F(x,y) w(y) \,dy
\end{equation}
holds.

\begin{thm}
When (\ref{D54}) through (\ref{D60}) hold, and when the integrals WRT $y$ in (\ref{D55}) and (\ref{D60}) are Riemann integrable, then
\begin{equation}
\label{D57}
\int^{\beta}_{\alpha} w(y) \left[\int^{\infty}_{a} \hspace{-19.9 bp} \textsf{Z} \hspace{10 bp} f(x,y) \,dx \right] dy = \int^{\infty}_{a} \hspace{-19.9 bp} \textsf{Z} \hspace{10 bp} \left[ \int^{\beta}_{\alpha} w(y) f(x,y) \,dy \right] \,dx .
\end{equation}
\end{thm}

In the proof below, we will again make extensive use of (\ref{D20}) and (\ref{D26}), and again simply use $c$ for the upper bound of all termination functions.

The restriction to Riemann integration over $y$ is due to the proof used below, but does not appear to be a fundamental limitation of our proposed integral definition.  A proof without this restriction is desirable, but unknown at this time.

\begin{proof}
The RHS of (\ref{D57}) becomes
\begin{equation}
\label{D61}
\int^{\infty}_{a} \hspace{-19.1 bp} \textsf{Z} \hspace{10 bp} \left[ g(x,\beta) - g(x,\alpha) \right] \,dx = \int^{\infty}_{a} \hspace{-19.1 bp} \textsf{Z} \hspace{10 bp} g(x,\beta) \,dx - \int^{\infty}_{a} \hspace{-19.1 bp} \textsf{Z} \hspace{10 bp} g(x,\alpha) \,dx
\end{equation}
\begin{eqnarray}
\label{D62}
&=& G(a,\alpha) + \lim_{b\to\infty}\left\{\int^{c}_{0} G(x+b,\alpha) \widetilde{z}_x(x,\alpha) \,dx\right\} \nonumber \\
& & - \left[G(a,\beta) + \lim_{b\to\infty}\left\{\int^{c}_{0} G(x+b,\beta) \widetilde{z}_x(x,\beta) \,dx\right\} \right]
\end{eqnarray}
\begin{eqnarray}
\label{D63}
&=& G(a,\alpha) + \lim_{b\to\infty}\left\{\int^{2c}_{0} G(x+b,\alpha) \widetilde{z}_x(x,\alpha) \otimes \widetilde{z}_x(x,\beta) \,dx\right\} \nonumber \\
& & - \left[G(a,\beta) + \lim_{b\to\infty}\left\{\int^{2c}_{0} G(x+b,\beta) \widetilde{z}_x(x,\alpha) \otimes \widetilde{z}_x(x,\beta) \,dx\right\} \right]
\end{eqnarray}
\begin{eqnarray}
\label{D64}
&=& \left[G(a,\alpha) - G(a,\beta)\right] + \nonumber \\
& & \lim_{b\to\infty}\left\{\int^{2c}_{0} \left[G(x+b,\alpha) - G(x+b,\beta)\right] \widetilde{z}_x(x,\alpha) \otimes \widetilde{z}_x(x,\beta) \,dx\right\}
\end{eqnarray}
\begin{eqnarray}
\label{D65}
&=& \left[G(a,\alpha) - G(a,\beta)\right] + \nonumber \\
& & \lim_{b\to\infty}\left\{\int^{2c}_{0} \left[G(x+b,\alpha) - G(x+b,\beta)\right] \right. \nonumber \\
& & \hspace{70 bp} \left. \widetilde{z}_x(x,\alpha) \otimes \widetilde{z}_x(x,\beta) \otimes \prod_{j = 1}^{M}{ \left\{ z_x(x,\tau _j) \right\} }\,dx\right\}
\end{eqnarray}
where the product symbol $\prod$ is taken as multiple convolutions, using (\ref{D26}) $M$ times.  Equation (\ref{D65}) holds for any arbitrary $M > 0$ and set $\left\{\tau _j\right\}$ subject to $\alpha \leq \tau _j \leq \beta$.  Using (\ref{D60}), (\ref{D65}) becomes
\begin{eqnarray}
\label{D66}
&=& \left[G(a,\alpha) - G(a,\beta)\right] + \nonumber \\
& & \lim_{b\to\infty}\left\{\int^{2c}_{0} \left[-\int^{\beta}_{\alpha} F(x+b,y) w(y) \,dy\right] \right. \nonumber \\
& & \hspace{70 bp} \left. \widetilde{z}_x(x,\alpha) \otimes \widetilde{z}_x(x,\beta) \otimes \prod_{j = 1}^{M}{ \left\{ z_x(x,\tau _j) \right\} }\,dx\right\} .
\end{eqnarray}

We now replace the integral using the standard Riemann integral limit definition.  That is, we will require that the interval $\alpha$ to $\beta$ be broken into $N$ intervals of length $\Delta y_i$, that each $y_i$ fall within the corresponding interval, and that the maximum length $\Delta y_i$ approach zero as $N$ approaches infinity,
\begin{eqnarray}
\label{D67}
&=& \left[G(a,\alpha) - G(a,\beta)\right] - \nonumber \\
& & \lim_{b\to\infty}\left\{\int^{2c}_{0} \lim_{N\to\infty}\left[\sum_{i=1}^{N} F(x+b,y_i) w(y_i) \,\Delta y_i\right] \right. \nonumber \\
& & \hspace{80 bp} \left. \widetilde{z}_x(x,\alpha) \otimes \widetilde{z}_x(x,\beta) \otimes \prod_{j = 1}^{M}{ \left\{ z_x(x,\tau _j) \right\} }\,dx\right\} .
\end{eqnarray}

For every $N$ and set of ${y_i}$, we can choose $M=N$ and $\tau _i = y_i$.  Thus, for every $N$ and $i$ in (\ref{D67}), we have that
\begin{eqnarray}
\label{D68}
& & \lim_{b\to\infty}\left\{\int^{2c}_{0} F(x+b,y_i) w(y_i) \,\Delta y_i \right. \nonumber \\
& & \hspace{80 bp} \left. \widetilde{z}_x(x,\alpha) \otimes \widetilde{z}_x(x,\beta) \otimes \prod_{j = 1}^{N}{ \left\{ z_x(x,y_j) \right\} }\,dx\right\}
\end{eqnarray}
exists.

We can then interchange the limit over $i$ and the integration in (\ref{D67}) to get
\begin{eqnarray}
\label{D69}
&=& \left[G(a,\alpha) - G(a,\beta)\right] - \nonumber \\
& & \lim_{b\to\infty}\left\{\lim_{N\to\infty}\left[\sum_{i=1}^{N} \int^{2c}_{0} F(x+b,y_i) w(y_i) \,\Delta y_i \right. \right. \nonumber \\
& & \hspace{70 bp} \left. \left. \widetilde{z}_x(x,\alpha) \otimes \widetilde{z}_x(x,\beta) \otimes \prod_{j = 1}^{N}{ \left\{ z_x(x,y_j) \right\} }\,dx \right]\right\} .
\end{eqnarray}

Using the existence of (\ref{D68}) again, we can interchange the two limits to get
\begin{eqnarray}
\label{D70}
&=& \left[G(a,\alpha) - G(a,\beta)\right] - \nonumber \\
& & \lim_{N\to\infty}\left[\lim_{b\to\infty}\left\{\sum_{i=1}^{N} \int^{2c}_{0} F(x+b,y_i) w(y_i) \,\Delta y_i \right. \right. \nonumber \\
& & \hspace{70 bp} \left. \left. \widetilde{z}_x(x,\alpha) \otimes \widetilde{z}_x(x,\beta) \otimes \prod_{j = 1}^{N}{ \left\{ z_x(x,y_j) \right\} }\,dx \right\}\right] .
\end{eqnarray}

Using (\ref{D09}), the LHS of (\ref{D57}) becomes
\begin{eqnarray}
\label{D71}
& & \hspace{-15 bp} \int^{\beta}_{\alpha} w(y) \left[-F(a,y) - \lim_{b\to\infty}\left\{\int^{c}_{0} F(x+b,y) z_x(x,y) \,dx\right\} \right] dy
\label{D72} \nonumber \\
& & \hspace{-15 bp} = -\int^{\beta}_{\alpha} w(y) F(a,y) \,dy -\int^{\beta}_{\alpha} \left[\lim_{b\to\infty}\left\{\int^{c}_{0} w(y) F(x+b,y) z_x(x,y) \,dx\right\} \right] \,dy \nonumber \\
& & \hspace{-15 bp} = \left[G(a,\alpha) - G(a,\beta)\right] \nonumber \\
& &\hspace{-10 bp}  -\int^{\beta}_{\alpha} \left[\lim_{b\to\infty}\left\{\int^{3c}_{0} w(y) F(x+b,y) z_x(x,y) \otimes \widetilde{z}_x(x,\alpha) \otimes \widetilde{z}_x(x,\beta) \,dx\right\} \right] \,dy . \nonumber \\
& & \hspace{ 1 bp}
\end{eqnarray}

Replacing the integration with the Riemann integral limit definition,
\begin{eqnarray}
\label{D73}
& & \hspace{-10 bp} =\left[G(a,\alpha) - G(a,\beta)\right] \nonumber \\
& & \hspace{-10 bp} -\lim_{N\to\infty}\left[\sum_{i=1}^{N} \left[\lim_{b\to\infty}\left\{\int^{3c}_{0} w(y_i) F(x+b,y_i) \right. \right. \right. \nonumber \\
& & \hspace{100 bp} \left. \left. \left. 
\begin{array}{l}
\!\!\!\! \\
\!\!\!\! \end{array}
\hspace{-20 bp} z_x(x,y_i) \otimes \widetilde{z}_x(x,\alpha) \otimes \widetilde{z}_x(x,\beta) \,dx \right\} \,\Delta y_i \right] \right] \hspace{ 1 bp} .
\end{eqnarray}

Combining the other $N-1$ additional termination functions with each $i$ using (\ref{D26}), we have
\begin{eqnarray}
\label{D74}
&=&\left[G(a,\alpha) - G(a,\beta)\right] \nonumber \\
& & -\lim_{N\to\infty}\left[\sum_{i=1}^{N} \left[\lim_{b\to\infty}\left\{\int^{3c}_{0} w(y_i) F(x+b,y_i) \right. \right. \right. \nonumber \\
& & \hspace{50 bp} \left. \left. \left. \widetilde{z}_x(x,\alpha) \otimes \widetilde{z}_x(x,\beta) \otimes \prod_{j = 1}^{N}{ \left\{ z_x(x,y_j) \right\} } \,dx\right\} \,\Delta y_i \right] \right] \hspace{ 1 bp}
\end{eqnarray}
or finally
\begin{eqnarray}
\label{D75}
 &=&\left[G(a,\alpha) - G(a,\beta)\right] \nonumber \\
& & -\lim_{N\to\infty}\left[\lim_{b\to\infty}\left\{\sum_{i=1}^{N} \left[\int^{3c}_{0} w(y_i) F(x+b,y_i) \right. \right. \right. \nonumber \\
& & \hspace{60 bp} \left. \left. \left. \widetilde{z}_x(x,\alpha) \otimes \widetilde{z}_x(x,\beta) \otimes \prod_{j = 1}^{N}{ \left\{ z_x(x,y_j) \right\} } \,dx  \,\Delta y_i \right] \right\} \right] .
\end{eqnarray}

Equations (\ref{D75}) and (\ref{D70}) are the same, so (\ref{D57}) holds.
\end{proof}

Since conventional integration is a special case of this definition, we immediately get
\begin{lem}
When the corresponding conventional and generalized integrations WRT $x$ in (\ref{D76}) exist, (\ref{D54}) through (\ref{D60}) hold, and the integrations WRT $y$ exist in the Riemann sense,
\begin{equation}
\label{D76}
\int^{\beta}_{\alpha} w(y) \left[\int^{\infty}_{a} \hspace{-19.9 bp} \textsf{Z} \hspace{10 bp} f(x,y) \,dx \right] dy = \int^{\infty}_{a} \left[ \int^{\beta}_{\alpha} w(y) f(x,y) \,dy \right] \,dx
\end{equation}
\end{lem}
and
\begin{lem}
When the corresponding conventional and generalized integrations WRT $x$ in (\ref{D77}) exist, (\ref{D54}) through (\ref{D60}) hold, and the integrations WRT $y$ exist in the Riemann sense,
\begin{equation}
\label{D77}
\int^{\beta}_{\alpha} w(y) \left[\int^{\infty}_{a} f(x,y) \,dx \right] dy = \int^{\infty}_{a} \hspace{-19.9 bp} \textsf{Z} \hspace{10 bp} \left[ \int^{\beta}_{\alpha} w(y) f(x,y) \,dy \right] \,dx
\end{equation}
\end{lem}

The existence of the conventional integral on the LHS of (\ref{D77}) may imply that using the generalized definition on the RHS is never necessary, but this has not been investigated.

An interesting open question is whether the existence of the LHS of (\ref{D57}) implies the existence of the RHS.  Intuition suggests that this is the case, but I have no proof.  It is clearly not true that the existence of the RHS of (\ref{D57}) implies the existence of the LHS.  Consider $f(x,y)$ constant, where $w(y)$ integrates to zero between $\alpha$ and $\beta$.

\section{Change of variable of integration}

There is a restriction on manipulating integrations under this definition that is not present in the conventional definition.  While some basic forms of integration by substitution can be used to change the variable of integration, arbitrary changes can not in general be performed.  To see this, begin with (\ref{D01}), rewritten here with $u$ as the variable of integration
\begin{equation}
\label{D78}
\int^{\infty}_{\alpha} \hspace{-19.1 bp} \textsf{Z} \hspace{10 bp} f(u) \,du = \lim_{\beta\rightarrow\infty}\left\{\int^{\beta}_{\alpha}f(u) \,du + \int^{\beta+\gamma}_{\beta} f(u) \zeta(u-\beta) \,du\right\} .
\end{equation}

Now let 
\begin{equation}
\label{D79}
u = u(x) \mbox{ } \mbox{ } \mbox{ } \mbox{ } x = u^{-1}(u)
\end{equation}
subject to the restrictions that $u(x)$ be a monotonic function of $x$, satisfying
\begin{equation}
\label{D80}
\lim_{x\rightarrow\infty}u(x) = \infty
\end{equation}
\begin{equation}
\label{D81}
\alpha = u(a) \mbox{ } \mbox{ } \mbox{ } \mbox{ } a = u^{-1}(\alpha)
\end{equation}
\begin{equation}
\label{D82}
\beta = u(b) \mbox{ } \mbox{ } \mbox{ } \mbox{ } b = u^{-1}(\beta)
\end{equation}

and let
\begin{equation}
\label{D83}
g(x) \equiv f(u(x)) \frac{du}{dx} .
\end{equation}

Substituting, the RHS of (\ref{D78}) becomes
\begin{equation}
\label{D84}
\lim_{\beta\rightarrow\infty}\left\{\int^{u^{-1}(\beta)}_{u^{-1}(\alpha)}f(u(x)) \frac{du}{dx} \,dx + \int^{u^{-1}(\beta + \gamma)}_{u^{-1}(\beta)} f(u(x)) \zeta(u(x)-u(b)) \frac{du}{dx} \,dx\right\} .
\end{equation}
Using (\ref{D80}) and (\ref{D82}), we can replace the limit $\beta\rightarrow\infty$ with the limit $b\rightarrow\infty$, and using (\ref{D83}), we get
\begin{equation}
\label{D85}
\lim_{b\rightarrow\infty}\left\{\int^{b}_{a}g(x) \,dx + \int^{b+c}_{b} g(x) z(x-b,b) \,dx\right\}
\end{equation}
where
\begin{equation}
\label{D86}
z(x-b,b) \equiv \zeta(u(x)-u(b)) .
\end{equation}

While (\ref{D85}) is equal to (\ref{D84}), and looks very similar to (\ref{D01}), we see that $z(\cdot)$ is an explicit function of $b$.  Variation of termination functions WRT $b$ is not allowed under our definition, since this can lead to non-unique values for the integrals.  To maintain equality in going from (\ref{D84}) to (\ref{D85}), an explicit variation WRT $b$ is required by (\ref{D86}).  The converse of this is that if a change of the variable of integration is performed, and a different termination function $z(x)$ that is \textit{not} a function of $b$ is used in (\ref{D85}), the value of the limit, and hence the integral, may be different.

In the following example, a change of variable using $u$-substitution is performed, affecting the value of the integral.

\begin{example}
\label{EX06}
 Consider a square wave function, with value alternating between $\pm 1$,
\[
f(x) \equiv 1 - 2 * \textrm{mod}\left( \left\lfloor x \right\rfloor , 2 \right)
\]
where $\left\lfloor x \right\rfloor$ is the floor function, and choose
\[
F(x) = \left\{ \begin{array}{lcc}
{ \mbox{ } \mbox{ } \left[ x - \left\lfloor x \right\rfloor \right] \mbox{ } \mbox{ }} {\mbox{ } \mbox{ } \mbox{ }} {\textrm{mod} \left( \left\lfloor x \right\rfloor , 2 \right) = 0} \\
{1 - \left[ x - \left\lfloor x \right\rfloor \right]} {\mbox{ } \mbox{ } \mbox{ }} {\textrm{mod}\left( \left\lfloor x \right\rfloor , 2 \right) = 1} \\
\end{array} \right. .
\]

Using (\ref{D09}),
\[
\textsf{Z} \hspace{-10.5 bp} \int_{0}^{\infty} f(x) \,dx = -F(0) - \lim_{b\to\infty}\left\{\int^{c}_{0} F(x+b) z'(x) \,dx\right\} .
\]
We can select $c = 1$ and use
\[
z(x) = 1/2 \mbox{ } \mbox{ } \mbox{ } 0 < x < 1
\]
so
\[
z'(x) = -1/2 \left( \delta(x) + \delta (x+1) \right)
\]
and we then get
\[
\textsf{Z} \hspace{-10.5 bp} \int_{0}^{\infty} f(x) \,dx = 1 / 2 .
\]

Now, attempt a change of variable.  Letting
\[
x = u + \alpha \sin(\pi u)
\]
\[
dx = \left( 1 + \alpha \pi \cos(\pi u) \right) du
\]
the integral becomes
\[
\textsf{Z} \hspace{-10.5 bp} \int_{0}^{\infty} f(u + \alpha \sin(\pi u)) \left( 1 + \alpha \pi \cos(\pi u) \right) \,du .
\]

If we restrict $\alpha$ to
\[
-\frac{1}{\pi} \leq \alpha \leq \frac{1}{\pi}
\]
then
\[
f(u + \alpha \sin(\pi u)) = f(u)
\]
and we get
\begin{eqnarray*}
& & \hspace{-20 bp} \textsf{Z} \hspace{-10.5 bp} \int_{0}^{\infty} f(u + \alpha \sin(\pi u)) \left( 1 + \alpha \pi \cos(\pi u) \right) \,du = \\
& & \textsf{Z} \hspace{-10.5 bp} \int_{0}^{\infty} f(u) \,du + \textsf{Z} \hspace{-10.5 bp} \int_{0}^{\infty} \alpha \pi f(u) \cos(\pi u)  \,du .
\end{eqnarray*}

Take
\[
g(u) \equiv \alpha \pi f(u) \cos(\pi u)
\]
and
\[
G(u) = \alpha \pi \sin \left( \pi \left( u - \left\lfloor u \right\rfloor \right) \right) .
\]

We can use $c=1$ and $z(u) = 1 - u$ giving 
\[
z'(u) = \left\{ \begin{array}{lcc}
{ -1 } {\mbox{ } \mbox{ } \mbox{ } \mbox{ } \mbox{ }} {0 < u < 1} \\
{ \mbox{ } 0 \mbox{ }} {\mbox{ } \mbox{ } \mbox{ } \mbox{ } \mbox{ }} {\textrm{else}} \\
\end{array} \right.\]
to find
\[
\textsf{Z} \hspace{-10.5 bp} \int_{0}^{\infty} g(u) \,du = 2 \alpha .
\]

We then get
\[
\textsf{Z} \hspace{-10.5 bp} \int_{0}^{\infty} f(u) \,du + \textsf{Z} \hspace{-10.5 bp} \int_{0}^{\infty} \alpha \pi f(u) \cos(\pi u)  \,du = 1/2 + 2 \alpha \neq \textsf{Z} \hspace{-10.5 bp} \int_{0}^{\infty} f(x) \,dx
\]
and using $u$-substitution to change the variable of integration can be seen not be valid.
\end{example}

In the restricted case where $u$ is linearly related to $x$, it can be shown that the resulting termination function is not a function of $b$, so a change of the variable of integration may be performed.  Assume
\begin{equation}
\label{D87}
u = r + sx \mbox{ } \mbox{ } \mbox{ } \mbox{ } x = \frac{u-r}{s}
\end{equation}
with $s > 0$.  Then
\begin{eqnarray}
\alpha = r+sa \mbox{ } \mbox{ } &\mbox{ }& \mbox{ } \mbox{ } a = (\alpha-r)/s \\
\beta = r+sb \mbox{ } \mbox{ } &\mbox{ }& \mbox{ } \mbox{ } b = (\beta-r)/s \\
\gamma = sc \mbox{ } \mbox{ } &\mbox{ }& \mbox{ } \mbox{ } c = \gamma/s .
\end{eqnarray}

Then (\ref{D84}) specializes to
\begin{eqnarray}
& &\lim_{\beta\rightarrow\infty}\left\{\int^{(\beta-r)/s}_{(\alpha-r)s}f(r+sx) s \,dx + \right. \nonumber \\
& & \hspace{70 bp} \left. \int^{((\beta + \gamma - r)/s)}_{(\beta-r)/s} f(r+sx) \zeta(r+sx-(r+sb)) s \,dx\right\}
\end{eqnarray}
or, since $b \rightarrow \infty$ as $\beta \rightarrow \infty$,
\begin{equation}
s \lim_{\beta\rightarrow\infty}\left\{\int^{b}_{a}f(r+sx) \,dx + \int^{b+c}_{b} f(r+sx) z(x-b) \,dx\right\}
\end{equation}
where
\begin{equation}
z(x-b) \equiv \zeta(sx-sb) .
\end{equation}

For this restricted case, the function $z(x)$ is not a function of $b$, and satisfies the criteria for a termination function.

\section{Conclusion}

The generalized definition of an improper integral with infinite bounds presented here is a more powerful alternative to the conventional definition.  The range of functions which are integrable to infinite limits is greatly expanded.

The generalized definition gives the same results as the conventional definition when that applies, and preserves uniqueness and linearity.  The new definition allows interchange of the order of differentiation and integration whenever the two integrals exist under the definition.  Also allowed is interchange of the order of integration of iterated integration, again whenever the integrations exist under the generalized definition.  The ability to rigorously interchange order of integrations, or order of integration and differentiation, in cases where integrals under the conventional definition do not converge, provides an added tool for manipulation of complicated integrals.  

The definition presented here is distinct from, but compatible with, the Hadamard finite part integral~\cite{Zwillinger:hadamard}.  It can be thought of as providing an intermediate level of convergence, less convergent than integrals which are absolutely or conditionally convergent, but more so than integrals which require finite part integration for removal of an infinite component.  Divergent integrands which can be handled using the Hadamard finite part definition are also divergent under our definition.  Conversely, integrands which require our definition are not integrable in the conventional Hadamard sense.  The two definitions can be used together if necessary, with our definition used for the evaluation of the finite part, after the infinite term is removed.

There are several open questions remaining for this definition.  It is unknown whether differentiation under the integral sign is always possible, or whether bringing another integration inside the improper integral is always possible.  It is also unknown whether the restriction to Riemann integrability for interchanging the order of integration can be removed.  Perhaps most importantly, the restrictions on termination functions may be stricter than necessary, and could be relaxed, further expanding the range of functions which are integrable under this definition.  A restriction on changing the variable of integration, not present under conventional integration, was demonstrated, although the precise restrictions that must be satisfied are unknown.  Answering these questions is desirable, and a subject for future research.  The results presented here, however, are independent of the answers.

If one desires to use a conventional definition of the improper integral for infinite limits, the development in this paper still allows more freedom to interchange order of integration and differentiation.  Since the generalized definition is defined in terms of standard integrals and limits, and since integrals using the generalized definition are equal to those using the conventional definition when the latter exists, one need only \emph{begin} with an integral that exists under the conventional definition.  Changes of order of integration and differentiation which later lead to integrands which would be nonconvergent may subsequently be performed, using the generalized definition as necessary.  All results obtained using the generalized definition remain equal to that of the original conventional definition.

\bibliographystyle{alpha}
\bibliography{ExtendedIntegrationFull}

\end{document}